\input amstex.tex
\input amsppt.sty   
\magnification 1200
\vsize = 8.5 true in
\hsize=6.2 true in
\NoRunningHeads 
\nologo       
\parskip=\medskipamount
        \lineskip=2pt\baselineskip=18pt\lineskiplimit=0pt
       
        \TagsOnRight
        \NoBlackBoxes

        \topmatter
        \title
        Spectral methods in PDE
        \endtitle
\author
         W.-M.~Wang        \endauthor        
\address
{D\'epartement de Math\'ematique, Universit\'e Paris Sud, 91405 Orsay Cedex, FRANCE}
\endaddress
        \email
{wei-min.wang\@math.u-psud.fr}
\endemail
\abstract
This is to review some recent progress in PDE. The emphasis is 
on (energy) supercritical nonlinear Schr\"odinger equations.
The methods are applicable to other nonlinear equations.   
 \endabstract
\endtopmatter
\document
\head{\bf 1. Introduction}\endhead
We consider the nonlinear Schr\"odinger equation on the $d$-torus $\Bbb T^d=[0, 2\pi)^d$:
$$
i\frac\partial{\partial t}u =-\Delta u+|u|^{2p}u+H(x,u,\bar u)\qquad (p\geq 1, p\in\Bbb N),\tag 1.1
$$
with periodic boundary conditions: $u(t,x)=u(t, x+2n\pi)$, $x\in [0, 2\pi)^d$ for all $n\in\Bbb Z^d$,
where $H(x,u,\bar u)$ is analytic and has the expansion:
$$H(x,u,\bar u)=\sum_{m=1}^\infty \alpha_m(x)|u|^{2p+2m} u,$$
with $\alpha_m$ periodic and uniformly real analytic. 
The integer $p$ in (1.1) is {\it arbitrary}. 

We study (1.1) from two perspectives, that of existence of quasi-periodic solutions and that 
of existence and uniqueness of smooth solutions to Cauchy problems. The link of these two questions are 
classical  as the known invariant measure for smooth solutions are supported on KAM tori. 
Conceptually, one could view quasi-periodic solutions as providing a  ``basis" for the nonlinear solutions.

It is well known \cite{Bo1} that on $\Bbb T^d$, when $H=0$, (1.1) is locally wellposed in $H^s$ for 
$$s>\max (0, \frac{1}{2}(d-\frac{2}{p})),\tag 1.2$$
which is derived by linearizing about the flow of the Laplacian and proving 
$L^p$ estimates of its eigenfunction solutions (Strichartz estimates). 

For $d\geq 3$ and sufficiently large $p$, the right side of (1.2) is larger than $1$. For example, in
dimension $4$, the quintic nonlinear Schr\"odinger equation is locally well-posed in $H^s$ for $s>1$, above the Hamiltonian 
$H^1$ topology, hence there is no available conservation law. These equations are therefore {\it supercritical} 
as there is no a priori global existence from patching up local solutions, not even for small data, as (1.1) is non-dispersive,
i. e., $\Vert u\Vert_\infty$ cannot tend to $0$ as $t\to\infty$ on the torus $\Bbb T^d$. 

The main purpose
of this review is to describe the new method developed in \cite {W2, 3} to construct global or almost global solutions to (1.1) both in the perturbative and semi-classical regime.  The semi-classical limit reveals further
the geometric nature of this construction. This is not surprising, as already in the critical case, the local existence
time depends on the Fourier geometry of the solution and not just its $H^1$ norm. 

Following the PDE custom, when considering applications to the Cauchy problem, we will set the
higher order terms $H=0$. This is also because for small data, the construction carries over verbatim to $H\neq 0$.
For KAM construction we choose to keep the translation invariance breaking term $H$, which traditionally lead to
more complicated small divisors in the normalizing transform. Our method is different, it hinges on establishing a spectral gap, which is  {\it stable} under small perturbations.  So we keep $H$ to illustrate this stability.

To understand this stability, it is important to remark that the spectral gap is in the ``space-time" sense as we put 
space and time on equal footing via a space-time Fourier series.

The point of departure here is that this spectral gap is produced by the nonlinearity itself and {\it not} by eigenvalue
variation of the linear operator, which has been the tradition. It is well known that eigenvalue variations
are difficult to achieve with a multiplicative potential in dimensions $2$ and above due to degeneracy 
of the Laplacian. 

A step was taken in \cite{W1} to break the degeneracy with a trigonometric polynomial potential
in dimension $2$, but much remains to be done.  On the other hand, the spectral gap 
produced by the nonlinearity is geometric in origin and hence robust. This is also 
because (1.1) essentially preserves translation invariance and hence the algebraic properties
of the exponentials, while adding a multiplicative potential destroys integrability.  

In the last part of the paper, we consider the nonlinear Schr\"odinger equation on general compact
manifolds. We obtain critical Sobolev exponents for local well-posedness on general surfaces. 
This is an indication that the results on the flat torus may not be improved geometrically.
\bigskip

\head{\bf 2. Quasi-periodic solutions}\endhead
In this section, we consider the KAM aspect of (1.1). The Kolmogorov non-degeneracy conditions (or its weaker
versions) are completely violated, as we perturb about a (infinite dimensional) linear system and {\it not} an integrable non-linear system.

Let $u^{(0)}$ be a solution to the linear equation: 
$$i\frac\partial{\partial t}u^{(0)} =-\Delta u^{(0)}.\tag 2.1$$
We seek quasi-periodic solutions to (1.1) with $b$ frequencies close to $u^{(0)}$ in the form of a nonlinear 
space-time Fourier series: 
$$
u(t, x)=\sum_{(n,j)}\hat u(n, j)e^{in\cdot\omega t}e^{ij\cdot x}, \qquad (n,j)\in\Bbb Z^{b+d},\qquad \tag 2.2
$$
with $\omega\in\Bbb R^b$ to be determined. 

Writing in this form, a solution $u^{(0)}$ to (2.1) 
with $b$ frequencies $\omega^{(0)}=\{j_k^2\}_{k=1}^{b}$ ($j_k\neq 0$) has Fourier support
$$\text{supp } {\hat u}^{(0)}=\{(-e_{j_k}, j_k), k=1,...,b\},$$
where $e_{j_k}$ is a unit vector in $\Bbb Z^b$ and $j_k\neq j_{k'}$ if $k\neq k'$. (Unless otherwise stated $j_k^2:=|j_k|^2$ etc.) 

We note that (2.2) treats space and time on equal footing, contrary to the usual ODE
approach, which views the solution as the time evolution of (Hamiltonian) vector fields. In this Hamiltonian
language, the vector fields corresponding to the nonlinear Schr\"odinger equation in (1.1) are infinite dimensional.   
We also note that in (2.2),  for each additional frequency in time, we add a dimension, which plays an important role
in untangling the resonances. 

Define the bi-characteristics
$$\Cal C=\{(n,j)\in\Bbb Z^{b+d}|\pm n\cdot\omega^{(0)}+j^2=0\}.\tag 2.3$$
$\Cal C$ is the solution set in the form (2.2) to (2.1) and its complex conjugate in the Fourier space.
This is the resonant set for the nonlinear equation (1.1).  The fact that 
the singularities are not isolated points is the major difficulty here. 

To treat this manifold of singularities, we observe that in the absence of the nonlinear term, there is no propagation
in the Fourier space $\Bbb Z^{b+d}$, as there is no coupling among different $(n, j)$, even though $\Cal C$
is an infinite set. In the Fourier space, the nonlinear term becomes a convolution operator and couples different
$(n, j)$. But as long as the connected (relative to this convolution) sets in $\Cal C$ remain bounded
in the first approximation, we should  still be able to construct quasi-periodic solutions close to a given linear solution
$u^{(0)}$ by implementing a Newton scheme. 

Let $$u^{(0)}(t, x)=\sum_{k=1}^b a_k e^{-ij_k^{2}t}e^{ij_k\cdot x},$$ be a solution to the linear equation (2.1). We achieve this by first making a {\it geometric} selection in the Fourier frequencies $\{j_k\}_{k=1}^b$ and then 
an amplitude selection in the Fourier coefficients $\{a_k\}_{k=1}^b$. The geometric selection is in order to 
ensure bounded connected sets in $\Cal C$ and is the main new ingredient in this theory.  

We call $u^{(0)}$ {\it generic}, if it satisfies the genericity conditions (i-iv) in \cite{W2}.   
Here it suffices to say that the genericity conditions pertain entirely to the Fourier support of $u^{(0)}$: $\{j_k\}_{k=1}^b\in(\Bbb R^d)^b $ and are determined by the $|u|^{2p}u$ term in (1.1) only. Moreover the non-generic set $\Omega$ is of codimension $1$ in $(\Bbb R^d)^b$.

The main result is 
\proclaim
{Theorem 1 \cite{W2}}
Assume $$u^{(0)}(t, x)=\sum_{k=1}^b a_k e^{-ij_k^{2}t}e^{ij_k\cdot x},$$ a solution to the linear equation (2.1) is generic and $a=\{a_k\}\in (0,\delta]^b=\Cal B(0, \delta)$. There exist $C$, $c>0$, such that for all
$\epsilon\in (0,1)$, there exists $\delta_0>0$ and for all $\delta\in(0,\delta_0)$ a Cantor set 
$\Cal G$ with 
$$\text{meas }\{\Cal G\cap \Cal B (0,\delta)\}/\delta^b\geq 1-C\epsilon^c.$$ 
For all $a\in\Cal G$, there is a quasi-periodic solution of $b$ frequencies to the nonlinear Schr\"odinger equation (1.1):
$$u(t, x)=\sum a_k e^{-i{\omega_k}t}e^{ij_k\cdot x}+\Cal O(\delta^{3}),$$
with basic frequencies $\omega=\{\omega_k\}$ satisfying 
$$\omega_k=j_k^2+\Cal O(\delta^{2p}).$$
The remainder $\Cal O(\delta^{3})$ is in an analytic norm about a strip of width $\Cal O(1)$ on $\Bbb T^{b+d}$.
\endproclaim

\noindent{\it Remark.} When $d=p=1$, the non-generic
set $\Omega=\emptyset$. {\it All} $u^{(0)}$ are generic and only amplitude selection is necessary.
This is the well understood scenario, see the Appendix in sect. 5.
To understand the substance of the geometric and amplitude excisions in the theorem, it is useful to take $H=0$ and note the perpetual existence of periodic solutions:
$$u=ae^{-i(j^{2}+|a|^{2p})t}e^{ij\cdot x}$$
to (1.1) for all $j\in\Bbb Z^d$ and $a\in\Bbb C$.

We have moreover the following semi-classical analog, which is new to the KAM context:
\proclaim{Corollary 1 \cite{W2}}
Set $H=0$ in (1.1). Assume $$u^{(0)}(t, x)=\sum_{k=1}^b a_k e^{-ij_k^{2}t}e^{ij_k\cdot x},$$ a solution to the linear equation (2.1) is generic, $\{j_k\}_{k=1}^b\in[K\Bbb Z^d]^b$, $K\in\Bbb N^+$  
and $a=\{a_k\}\in (0,1]^b=\Cal B(0, 1)$. There exist $C$, $c>0$, such that for all
$\epsilon\in (0,1)$, there exists $K_0>0$ and for all $K>K_0$ a Cantor set 
$\Cal G$ with 
$$\text{meas }\{\Cal G\cap \Cal B (0,1)\}\geq 1-C\epsilon^c.$$ 
For all $a\in\Cal G$, there is a quasi-periodic solution of $b$ frequencies to the nonlinear Schr\"odinger equation (1.1): 
$$u(t, x)=\sum a_k e^{-i{\omega_k}t}e^{ij_k\cdot x}+\Cal O(1/K^2),$$
with basic frequencies $\omega=\{\omega_k\}$ satisfying 
$$\omega_k=j_k^2+\Cal O(1).$$
The remainder $\Cal O(1/K^2)$ is in an analytic norm about a strip of width $\Cal O(1)$ in $t$ and 
$\Cal O(1/K)$ in $x$ on $\Bbb T^{b+d}$.
\endproclaim

\noindent{\it Remark.} These are quantitative, global,  $\Bbb L^2$ size $1$ and large (kinetic) energy solutions, which could be relevant to the compressible Euler equations. 

\bigskip
\noindent{\it A sketch of proof of the theorem}

We write (1.1) in the Fourier space, it becomes
$$\text{diag }(n\cdot\omega+j^2)\hat u+(\hat u*\hat v)^{*p}* \hat u+ \sum_{m=1}^\infty\hat\alpha_m*(\hat u*\hat v)^{*(p+m)}* \hat u=0,\tag 2.4$$
where $(n,j)\in\Bbb Z^{b+d}$, $\hat v=\hat{\bar u}$, $\omega\in\Bbb R^b$ is to be determined and
$$|\hat\alpha_m(\ell)|\leq C'e^{-c'|\ell|}\quad (C', c'>0)$$
for all $m$. 
From now on we work with (2.4), for simplicity we drop the hat and write $u$ for $\hat u$ and $v$ for $\hat v$ etc.
We seek solutions close to the linear solution $u^{(0)}$ of $b$ frequencies, 
$\text{supp } {u}^{(0)}=\{(-e_{j_k}, j_k), k=1,...,b\},$ with frequencies
$\omega^{(0)}=\{j_k^2\}_{k=1}^{b}$ ($j_k\neq 0$) 
and small amplitudes $a=\{a_k\}_{k=1}^b$ satisfying $\Vert a\Vert=\Cal O(\delta)\ll 1$.

We complete (2.4) by writing the equation for the complex conjugate. So we have 
$$
\cases
\text{diag }(n\cdot\omega+j^2)u+(u*v)^{*p}* u+\sum_{m=1}^\infty\alpha_m*(u*v)^{*(p+m)}*u=0,\\
\text{diag }(-n\cdot\omega+j^2)v+(u*v)^{*p}* v+\sum_{m=1}^\infty\alpha_m*(u*v)^{*(p+m)}*v=0,
\endcases\tag 2.5
$$
By supp, we will always mean the Fourier support, so we write $\text{supp } u^{(0)}$ for 
$\text{supp } {\hat u}^{(0)}$ etc. Let 
$$\Cal S=\text{supp } u^{(0)}\cup\text{supp } {\bar u}^{(0)}.$$

Denote the left side of (2.5) by $F(u, v)$. We make a Lyapunov-Schmidt decomposition into the $P$-equations:
$$ F(u, v)|_{\Bbb Z^{b+d}\backslash\Cal S}=0,$$
and the $Q$-equations:
$$ F(u, v)|_{\Cal S}=0.$$
We seek solutions such that 
$u|_\Cal S=u^{(0)}$. 
The $P$-equations are infinite dimensional and determine $u$ in the complement of $\text{supp }u^{(0)}$; 
the $Q$-equations are $2b$ dimensional and determine the frequency $\omega=\{\omega_k\}_{k=1}^b$. 

This Lyapunov-Schmidt method was introduced by Craig and Wayne \cite{CW} to construct periodic solutions
for the wave equation in one dimension. It was inspired by the multiscale analysis of Fr\"ohlich and Spencer \cite{FS}.  
The construction was further developed by Bourgain to embrace the full generality of quasi-periodic solutions
and in arbitrary dimensions $d$ \cite{Bo4, 6}. More recently, Eliasson and Kuksin \cite{EK} developed 
a KAM theory in the Schr\"odinger context. All the above results, however, pertain to parameter dependent tangentially non-resonant equations. 

We use a Newton scheme to solve the $P$-equations, with $u^{(0)}$ as the initial approximation. The major
difference with \cite{Bo4, 6, CW, EK} is that (2.5) is completely resonant and there are {\it no} parameters at this initial
stage. The frequency $\omega^{(0)}$ is an {\it integer} in $\Bbb Z^b$. So we need to proceed differently.

First recall the formal scheme: the first correction
$$\Delta \pmatrix u^{(1)}\\v^{(1)}\endpmatrix= \pmatrix u^{(1)}\\v^{(1)}\endpmatrix- \pmatrix u^{(0)}\\v^{(0)}\endpmatrix
=[F'( u^{(0)}, v^{(0)}]^{-1} F( u^{(0)}, v^{(0)}),\tag 2.6$$
where $\pmatrix u^{(1)}\\v^{(1)}\endpmatrix$ is the next approximation and $F'( u^{(0)}, v^{(0)})$ is the linearized 
operator on $\ell^2(\Bbb Z^{b+d}\times \Bbb Z^{b+d})$
$$F'=D+A,$$
where 
$$
D =\pmatrix \text {diag }(n\cdot\omega+j^2)&0\\ 0& \text {diag }(-n\cdot\omega+j^2)\endpmatrix$$
and
$$\aligned
A&=\pmatrix (p+1)(u*v)^{*p}& p(u*v)^{*p-1}*u*u\\ p(u*v)^{*p-1}*v*v& (p+1)(u*v)^{*p}\endpmatrix+\Cal O(\delta^{2p+2})
\quad  (p\geq 1),\\
&=A_0+\Cal O(\delta^{2p+2}).\endaligned$$
with $\omega=\omega^{(0)}$, $u=u^{(0)}$ and $v=v^{(0)}$.

Since we look at small data, $\Vert A\Vert=\Cal O(\delta^{2p})\ll 1$ and the diagonal: $\pm n\cdot\omega+j^2$
are integer valued, using the Schur complement reduction \cite{S1, 2}, the spectrum of $F'$ around $0$ is equivalent to that of a reduced operator on $\ell^2(\Cal C)$, where $\Cal C$ is defined in (2.3) and to $\Cal O(\delta^{2p+2})$ it is 
the same as the spectrum of $A_0$ on $\ell^2(\Cal C)$.

To implement the Newton scheme using (2.5), we need to bound $A_0^{-1}$, which leads to genericity conditions (i-iv) in \cite{W2}. From previous considerations,
it suffices to consider $A_0$ restricted to $\Cal C$. For generic $u^{(0)}$,  
$A_0|_{\Cal C}=\oplus\Cal A_0$, where $\Cal A_0$ are T\"oplitz matrices of sizes at most $(2b+d)\times (2b+d)$. 
This can be seen as follows.

Let 
$$\aligned &\Cal C^+=\{(n,j) |n\cdot\omega^{(0)}+j^2=0, j\neq 0\}\cup \{(n,0) |n\cdot\omega^{(0)}=0,n_1\leq 0\}, \\
&\Cal C^-=\{(n,j) |-n\cdot\omega^{(0)}+j^2=0, j\neq 0\}\cup \{(n,0) |n\cdot\omega^{(0)}=0,n_1>0\}, \\
&\Cal C^+\cap\Cal C^-=\emptyset ,\quad \Cal C^+\cup\Cal C^-=\Cal C.\endaligned\tag 2.7$$

Assume $(n,j)\in\Cal C^+$ is connected to $(n',j')\in\Cal C$ by the convolution operator $A_0$, then
$n'=n+\Delta n$ and $j'=j+\Delta j$, where $(\Delta n, \Delta j)\in\text{supp } (u^{(0)}*v^{(0)})^{*p}$,
if $(n',j')\in\Cal C^+$ and
$$
\cases
(n\cdot\omega^{(0)}+j^2)=0,\\
(n+\Delta n)\cdot\omega^{(0)}+(j+\Delta j)^2=0;
\endcases\tag 2.8 
$$
and if $(n',j')\in\Cal C^-$, then $(\Delta n, \Delta j)\in\text{supp } (u^{(0)}*v^{(0)})^{*p-1}*u^{(0)}*u^{(0)}$ and
$$
\cases
(n\cdot\omega^{(0)}+j^2)=0,\\
-(n+\Delta n)\cdot\omega^{(0)}+(j+\Delta j)^2=0.
\endcases\tag 2.9
$$
(Clearly the situation is similar if  $(n,j)\in\Cal C^-$.)

(2.8, 2.9) define a system of polynomial equations. The generic condition on $u^{(0)}$ is imposed precisely to ensure that, aside from the degenerate case, (2.8, 2.9) have a solution for at most $d+2$ equations. The degenerate case can be analyzed and correspond to systems
of at most $2b$ equations reflecting translation invariance. 
So the largest connected set is of size at most $\text {max }(2b, d+2)\leq 2b+d$. Analysis of the degenerate case is 
in fact the reason for requiring the leading nonlinear $\Cal O(\delta^{2p+1})$ term in (1.1) to be independent of $x$. The $x$ dependence of the higher order terms do not matter as they are treated as perturbations.

The invertibility of $A_0$ is then ensured by making an initial excision in $a$ as $0$ is typically {\it not}
an eigenvalue of a finite matrix. So 
$\Vert {F'}^{-1}\Vert\asymp\Vert A_0^{-1}\Vert \leq\Cal O (\delta^{-2p})$. Let 
$$F_0(u^{(0)}, v^{(0)})= \pmatrix (u^{(0)}*v^{(0)})^{*p}*u^{(0)}\\(u^{(0)}*v^{(0)})^{*p}*v^{(0)}\endpmatrix.$$ 
 By requiring 
 $$\text{supp }F_0(u^{(0)}, v^{(0)})\cap\{\Cal C\backslash\Cal S\}=\emptyset,$$
 which is part of the genericity conditions in \cite{W2}, we obtain from (2.6)
 $$\Vert \Delta u^{(1)}\Vert=\Vert \Delta v^{(1)}\Vert\leq \Cal O(\delta^{3})$$
for small $\delta$. Inserting this into the $Q$-equations, which determine $\omega$, we achieve amplitude-frequency
modulation:
$$\aligned \Vert\Delta\omega^{(1)}\Vert &\asymp\Cal O(\delta^{2p})\\
\big| \det(\frac{\partial \omega^{(1)}}{\partial a})\big|&\asymp \Cal O(\delta^{2p-1})>0\endaligned$$
ensuring transversality and moreover Diophantine $\omega^{(1)}$ on a set of $a$ of positive measure.
The tangentially non-resonant perturbation theory in \cite {Bo4, 6} becomes available.

Previously, quasi-periodic solutions were constructed using partial Birkhoff normal forms for the resonant Schr\"odinger equation in the presence of the cubic nonlinearity in dimensions $1$ and $2$ \cite{GXY, KP}. These constructions rely on the specifics of the resonance geometry given by the cubic nonlinearity, see the Appendix in sect. 5. 
\bigskip

\head{\bf 3. Almost global existence}\endhead
In this section, we set $H=0$. As in sect. 2, we have results both on the perturbative as well as the semi-classical case.  For the perturbative case, it is convenient to add a small parameter $\delta$ and look at data of size $1$. 
So we investigate the Cauchy problem on $\Bbb T^d$:
$$
\cases i\frac\partial{\partial t}u =-\Delta u+\delta |u|^{2p}u\qquad (p\geq 1, p\in\Bbb N \text{ \it arbitrary}),\\
u(t=0)=u_0,\endcases
\tag 3.1
$$

Relying on the geometric information afforded by the resonance analysis in sect. 2 and linearizing
about a suitable approximate quasi-periodic solutions, we prove the following:
\proclaim{Theorem 2}
Let $u_0=u_1+u_2$. Assume $u_1$ is {\it generic} satisfying (I. i-iv.) and $\Vert u_2\Vert=\Cal O(\delta)$, 
where $\Vert\cdot\Vert$ is an analytic norm (about a strip of width $\Cal O(1)$) on $\Bbb T^d$. Let $\Cal B (0,1)=(0, 1]^b$,
where $b$ is the dimension of the Fourier support of $u_1$. Then there exists an open set $\Cal A\subset
\Cal B (0,1)$ of positive measure, such that for all $A>1$, there exists $\delta_0>0$, such that for 
all $\delta\in(-\delta_0, \delta_0)$, if $\{|\hat u_1|\}\in\Cal A$, then (1.1) has a unique solution $u(t)$ for 
$|t|\leq\delta^{-A}$ satisfying $u(t=0)=u_0$ and $\Vert u(t)\Vert  \leq \Vert u_0\Vert+\Cal O(\delta)$, moreover
meas $\Cal A\to 1$ as $\delta\to 0$.   
\endproclaim

\noindent{\it Remark.} It is essential that the set $\Cal A$ is {\it open}, which enables us to establish an open mapping theorem to analyze Cauchy problems .

For perturbations of the $1d$ cubic NLS ($d=p=1$), similar stability results are proven in \cite{Ba, Bo5}. For parameter
dependent equations see \cite{BG, Bo3}. The equations treated in \cite{Ba, BG, Bo3, 5} are either $\Bbb L^2$ 
or  essentially $\Bbb L^2$ well-posed. So there is a priori global existence. 

The equations treated in the theorem are of a different nature, there is no a priori global existence
from conservation laws. In fact existence is obtained via explicit construction. Linearizing about approximate quasi-periodic solutions to prove existence and uniqueness for a 
time arbitrarily longer than local existence time, which is $\Cal O(\delta^{-1})$, is the main novelty. 

As in sect. 2, we also have the following semi-classical counterpart, providing quantitative, almost global, $\Bbb L^2$ size $1$
and large (kinetic) energy solutions to Cauchy problems. These solutions could be relevant to Cauchy problems
for compressible Euler equations.

\proclaim{Corollary 2}
Set $\delta=1$ in (3.1). Assume $u_0$ is {\it generic} with frequencies $\{j_k\}_{k=1}^b\in [K\Bbb Z^d]^b$, $K\in\Bbb N^+$.
Let $\Cal B (0,1)=(0, 1]^b$. Then there exists an open set $\Cal A\subset
\Cal B (0,1)$ of positive measure, such that for all $A>1$, there exists $K_0>0$, such that for 
all $K>K_0$, if $\{|\hat u_0|\}\in\Cal A$, then (1.1) has a unique solution $u(t)$ for 
$|t|\leq K^{A}$ satisfying $u(t=0)=u_0$ and $\Vert u(t)\Vert  \leq \Vert u_0\Vert+\Cal O(1/K^2)$, 
where $\Vert\cdot\Vert$ is an analytic norm (about a strip of width $\Cal O(1/K)$) on $\Bbb T^d$, 
moreover meas $\Cal A\to 1$ as $K\to \infty$.   
\endproclaim

\noindent{\it Remark.} The previous related results are up to time  $\Cal O(K)$, by solving the associated Hamilton-Jacobi equations before the arrival of caustics, cf. \cite{Car}.

\bigskip
\noindent{\it A sketch of the proof}

Writing the first equation in (3.1) as $F(u)=0$, for $u_0$ satisfying the conditions in the theorem,
we first find an approximate solution $v$ such that 
$$\cases F(v)=\Cal O(\delta^r),  \qquad\qquad\qquad\qquad(3.2)\\
v(t=0)-u_0=\Cal O(\delta^r), \,\quad\quad\qquad\,(3.3)\endcases$$
where $r>A>1$.

This approximate solution $v$ is quasi-periodic with $\Cal O(|\log\delta|)$ number of basic frequencies.
The construction of $v$ comprises of a {\it finitely} iterated Newton scheme to construct quasi-periodic solutions
which solve (3.2) but (3.3) only to order $\delta$ and then establishing an open mapping theorem using the spectral gap to solve (3.3) to order $\delta^r$.

Differentiate these quasi-periodic solutions with respect to the Fourier coefficients of $u^{(0)}$ gives 
a basis which spans $\Bbb L^2(\Bbb T^d)$ and allows to control the linearized flow.
Schematically this could be understood as follows. 

Assume $u$ is a solution satisfying the equation 
$F(u)=0$ and that it depends on a parameter $a$, then $\partial u/\partial a$ is a solution 
to the linearized equation: $$F'(u)(\frac{\partial u}{\partial a})=0.$$
The difficulty here is that in the linearized equation for $u$, there is $\bar u$. However the
generic condition allows us to control the coupling of $u$ and $\bar u$ on the flat torus.  

Using Duhamel's formula and the linearized flow to control
the difference of (3.1) and (3.2, 3.3), we conclude the proof of the theorem.
 \bigskip
 
 \head{\bf 4. Nonlinear Schr\"odinger equations on surfaces}\endhead
 We now consider the problem of nonlinear Schr\"odinger equations on smooth compact manifolds (without boundary):
 $$
\cases i\frac\partial{\partial t}u =-\Delta u+|u|^{2p}u,\\
u(t=0)=u_0.\endcases
$$
We specialize to dimension $2$ and cubic nonlinearity ($p=1$) as this is the most amenable. So we study
the Cauchy problem
$$
\cases i\frac\partial{\partial t}u =-\Delta u+|u|^{2}u,\\
u(t=0)=u_0, \endcases\tag 4.1
$$
on a compact $2$-manifold $M$.   

It is known from \cite{BGT1, ST} that (4.1) is uniformly locally wellposed in $H^s$ for $s>1/2$, in fact the flow 
is Lipshitz. We note that the corresponding Strichartz estimate on the flat $2$-torus gives $s>0$, cf. (1.2) and
on the $2$-sphere $s>1/4$ \cite{BGT2}. Both are known to be optimal \cite {CCT, BGT2}. 
However the optimality of $1/2$ on general surfaces $M$ remained open. 

The following proposition constructs a counter example on the torus of revolution:
 \proclaim {Proposition}
 Let $ds^2=dx^2+g(x)dy^2$ with $g\in C^3$ and admitting a unique global maximum.  Then there
 are initial data, which are eigenfunctions of the Laplacian, such that the flow map is not Lipshitz
 in $H^s$ for $s<1/2$. 
 \endproclaim

So the torus of revolution provides the geometric obstruction for going below $H^{1/2}$ and $1/2$ is the critical Sobolev exponents for local well-posedness on general surfaces. 

\noindent{\it A sketch of the proof}

The Laplace-Beltrami operator decomposes into a direct sum of one dimensional Schr\"odinger operators
indexed by the Fourier variable $k$ in the $y$ direction. We take as initial data the ground state $\psi$  of
such a Schr\"odinger operator when $k$ is large. Let $\lambda$ be its eigenvalue, the proof hinges on the 
equivalence estimate on the $L^\infty$ norm of $\psi$, namely
$$\Vert\psi\Vert_\infty\asymp\lambda^{1/8}\tag 4.2$$
 
The proof of (4.2) uses Fourier and Hermite series analysis and a lemma of Bourgain \cite{Bo2, D}, which relates $L^2$
and $L^\infty$ estimates. Afterwards one uses standard short time nonlinear analysis to complete the
proof. (For the analysis details see \cite{Cat}.)  
 
We note that the above instability is in the sense of derivatives by exploiting the fact that in the
linearized equation for $u$, there is $\bar u$, and there is growth of $L^\infty$ norm. On the flat
torus, the exponentials have $L^p$ norm $1$ for all $p$, which prevent this instablity.

The above instability is uniquely infinite dimensional as it relies on non-equivalence of norms. However  
 it is rooted in classical mechanics, namely the existence of lower dimensional tori, which leads to
 eigenfunction concentration phenomenon. The Hermite series precisely represents the missing 
 dimension in this picture. 
 
\noindent{\it Remark.} Using the same type of construction, $2/3$ \cite {BSS} should be optimal for surfaces
with boundary.  
 
\head{\bf 5. Appendix: the cubic nonlinearity on $\Bbb T^d$}\endhead 
For simplicity we write $u$ for $u^{(0)}$ and $\omega$ for $\omega^{(0)}$, the solutions and frequencies
of the linear equation. The symbols of convolution for the cubic nonlinearity are $|u|^2$, $u^2$ and ${\bar u}^2$.
Assume $(n,j)\in\Cal C^+$ ($(n,j)\in\Cal C^-$ works similarly). In order that $(n,j)$ is connected 
to $(n'j')\in\Cal C$, it is necessary that either
\item{(a)} $[u*v](n,j;n'j')\neq 0$ or
\item{(b)} $[u*u](n,j;n'j')\neq 0$.

\noindent Case (a): Since $$\aligned n\cdot\omega+j^2&=0,\\
n'\cdot\omega+{j'}^2&=0,\endaligned$$
subtracting the two equations gives immediately
$$(j_k-j_k')\cdot(j+j_k)=0,\tag 5.1$$
where $j_k, j_{k'}\in\Bbb Z^d$ ($k, k'=1,..., b$) and $j_k\neq j_{k'}$ if $k\neq k'$, are the $b$ Fourier components of $u$.

\noindent Case (b): Since $$\aligned n\cdot\omega+j^2&=0,\\
-n'\cdot\omega+{j'}^2&=0,\endaligned$$
adding the two equations gives immediately
$$(j+j_k)\cdot(j+j_{k'})=0,\tag 5.2$$
where $j_k, j_{k'}\in\Bbb Z^d$ ($k, k'=1,..., b$) and $j_k\neq j_{k'}$ if $k\neq k'$, are the $b$ Fourier components of $u$.

(5.1, 5.2) are precisely the well known resonant set for the partial Birkhoff normal form
transform in \cite{Bo4, GXY, KP}. (5.1, 5.2) describe rectangular type of geometry. 
$$\text{supp }F_0(u, v)\cap\{\Cal C\backslash\Cal S\}=\emptyset$$
for the cubic nonlinearity in any $d$. 
When $d=1$, (5.1, 5.2) reduce to a finite set of $2b$ lattice points in $\Bbb Z$: $\{j=\pm j_k$, $k=1,...,b\}$ 
and $\Omega=\emptyset$ in the theorems and corollaries.

\enddocument
\end